\documentclass[12pt]{article}
\usepackage{amssymb,amsmath,amsfonts,amsthm,amscd,latexsym,verbatim,graphics,epsfig,indentfirst,geometry}
\def\id{\mathrm{id}}
\def\pre{\mathrm{pre}}
\def\post{\mathrm{post}}
\def\RB{\mathrm{RB}}
\def\Var{\mathrm{Var}}
\def\Com{\mathrm{Com}}

\def\shuffle{\sqcup\mathchoice{\mkern-7mu}{\mkern-7mu}{\mkern-3.2mu}{\mkern-3.8mu}\sqcup}

\begin{document}

\begin{center}
{\Large
Universal enveloping commutative Rota---Baxter algebras of
precommutative and postcommutative algebras}

V. Gubarev
\end{center}

\begin{abstract}
Universal enveloping commutative Rota---Baxter algebras of
pre- and postcom\-mu\-ta\-tive algebras
are constructed. The pair of varieties $(\RB_\lambda\Com,\post\Com)$
is proved to be a PBW-pair and the pair $(\RB\Com,\pre\Com)$ is not.

\medskip
{\it Keywords}: Rota---Baxter algebra,
universal enveloping algebra, PBW-pair of varieties,
Zinbiel algebra, precommutative algebra, postcommutative algebra.
\end{abstract}

\section*{Introduction}

Linear operator $R$ defined on an algebra $A$ over the key field $\Bbbk$
is called Rota---Baxter operator (shortly RB-operator) of a weight $\lambda\in\Bbbk$
if it satisfies the relation
\begin{equation}\label{eq:RB}
R(x)R(y) = R( R(x)y + xR(y) + \lambda xy), \quad x,y\in A.
\end{equation}
An algebra with given RB-operator acting on it
is called Rota---Baxter algebra (RB-algebra).

G. Baxter defined commutative RB-algebra in 1960 \cite{Baxter60},
the relation \eqref{eq:RB} with $\lambda = 0$
is simply a generalization of by part integration formula.
J.-C. Rota, P. Cartier and others studied \cite{Rota68,Cartier72,Miller}
combinatorial properties of RB-operators and RB-algebras.
In 1980s, the deep connection between Lie RB-algebras and
classical Yang---Baxter equation was found \cite{BelaDrin82,Semenov83}.
In 2000, M. Aguiar showed \cite{Aguiar00} that a solution of associative
Yang---Baxter equation \cite{Zhelyabin} gives rise to a structure of associative RB-algebra.
There are a lot of applications of RB-operators in mathematical physics,
combinatorics, number theory, operads \cite{Kreimer,QFT00,NumThe,GubKol2013}.

There exist many constructions of free commutative RB-algebra
including ones given by J.-C. Rota, P. Cartier,
and L. Guo jointly with W. Keigher \cite{Rota68,Cartier72,GuoKeigher}.
Remark that the last one \cite{GuoKeigher} could be used to define
a structure of Hopf algebra on free commutative RB-algebra \cite{Hopf}.
In 2008, K. Ebrahimi-Fard and L. Guo obtained free associative RB-algebra \cite{FardGuo08}.
Different linear bases of free Lie RB-algebra were recently found \cite{Gub2016,GubKol2016,Chen16}.

In 1995 \cite{Loday95}, J.-L. Loday introduced algebras which satisfy an identity
$$
(x_1\succ x_2 + x_2\succ x_1)\succ x_3 = x_1\succ (x_2\succ x_3).
$$
We will call such algebras as precommutative algebras because
they play an analogous role for so called prealgebras as play commutative algebras
for ordinary algebras. In literature, they are also known as
dual Leibniz algebras (by Koszul duality) or Zinbiel algebras
(the word ``Leibniz'' written in the inverse order).
About precommutative algebras see, e.g. \cite{Nilpotency,Umirbaev,Zinbiel}.

In 1999 \cite{Dialg99}, J.-L. Loday introduced dendriform algebra
(we will call it as preassociative algebra).
A linear space endowed with two bilinear products $\succ,\prec$
is called a preassociative algebra if the following identities are satisfied:
$$
\begin{gathered}
(x_1\succ x_2+x_1\prec x_2)\succ x_3 = x_1\succ (x_2 \succ x_3), \quad
(x_1\succ x_2)\prec x_3=x_1\succ(x_2\prec x_3), \\
x_1\prec(x_2\succ x_3+x_2\prec x_3)=(x_1\prec x_2)\prec x_3.
\end{gathered}
$$
Given a preassociative algebra $A$, if we have $x\succ y = y\prec x$
for any $x,y\in A$, then $A$ is a precommutative algebra due to the product $\succ$.
The same space $A$ under the product $x\cdot y = x\succ y - y\prec x$
is a pre-Lie algebra \cite{Vinberg63,Gerst63,Koszul61}, i.e.
satisfies the identity
$(x_1 x_2)x_3 - x_1(x_2 x_3) = (x_2 x_1) x_3 - x_2 (x_1 x_3)$.

In \cite{Trialg01,Loday2007,Vallette2007},
postassociative, postcommutative, and postLie algebras were introduced.
All of them have an additional product and satisfy some certain identities.

The common definition of the varieties of pre- and post-$\Var$-algebra
for a variety $\Var$ could be found in \cite{BBGN2012,GubKol2014}.

In 2000 \cite{Aguiar00}, M. Aguiar noticed that
any commutative algebra with defined on it
a Rota---Baxter operator $R$ of zero weight
with respect to the operation $a\succ b = R(a)b$
is a precommutative algebra.
In 2007 \cite{Loday2007}, J.-L. Loday stated that
a commutative algebra $\langle A,\cdot\rangle$ with
an RB-operator of unit weight
is a postcommutative algebra under the operations
$\succ$ and $\cdot$ where $x\succ y = R(x)\cdot y$.

In 2013 \cite{BBGN2012}, this connection between RB-algebras and pre- and
postalgebras was generalized to any variety.
In 2013 \cite{GubKol2013}, it was proved that
any pre-$\Var$-algebra (post-$\Var$-algebra)
injectively embeds into its universal enveloping
RB-algebra of variety $\Var$ and (non)zero weight.

Based on the last result, we have a problem:
To construct an universal enveloping RB-algebra of a variety $\Var$
for pre- and post-$\Var$-algebra.
Another close problem is the following:
Whether the pairs of varieties $(\RB\Var,\mathrm{pre}\Var)$
and $(\RB_\lambda\Var,\mathrm{post}\Var)$, $\lambda\neq0$, are PBW-pairs \cite{PBW}?
Here by $\RB\Var$ ($\RB_\lambda\Var$) we mean the variety of
RB-algebras of a variety $\Var$ and (non)zero weight $\lambda$.

The questions stated above appeared in associative case
in the comments to the chapter V of the unique monograph
on RB-algebras by L. Guo \cite{Guo2011}
and was solved by the author in \cite{Gub2017}.
Also, see a brief history of the subject in associative
case in \cite{Gub2017}.

The current work is devoted to the solution of the stated problems
in commutative case. Although, the main method of the solution
in associative and commutative cases is similar,
the technical tools and the constructions are rather different
and it is hard to derive any of them from another one.

In \S1, we show connections between RB-algebras and classical,
modified and associative versions of Yang---Baxter equation,
we give preliminaries on pre- and postcommutative algebras, PBW-pairs of varieties.
Universal enveloping RB-algebras of precommutative (\S2)
and postcommutative (\S3) are constructed.
As corollaries we state that the pair of varieties
$(\RB_\lambda\Com,\mathrm{post}\Com)$ is a PBW-pair and
$(\RB\Com,\mathrm{pre}\Com)$ is not.

\section{Preliminaries}

\subsection{RB-operator}
Let us consider some well-known examples of RB-operators (see, e.g., \cite{Guo2011}):

{\bf Example 1}.
Given an algebra $A$ of continuous functions on $\mathbb{R}$,
an integration operator
$R(f)(x) = \int\limits_{0}^x f(t)\,dt$
is an RB-operator on $A$ of zero weight.

{\bf Example 2}.
Given an invertible derivation $d$ on an algebra $A$,
$d^{-1}$ is an RB-operator on $A$ of zero weight.

{\bf Example 3}.
Let $A = \{(a_1,a_2,\ldots,a_k,\ldots)\mid a_i\in \Bbbk\}$ be a infinite sum of a field $\Bbbk$.
An operator $R$ defined as
$R(a_1,a_2,\ldots,a_k,\ldots) = (a_1,a_1+a_2,\ldots,\sum\limits_{i=1}^k a_i,\ldots)$
is an RB-operator on $A$ of weight $-1$.

Further, if nothing else is specified,
by RB-operator we will mean an RB-operator of zero weight.

\subsection{Yang---Baxter equation}

Let $\mathfrak{g}$ be a semisimple finite-dimensional Lie algebra over $\mathbb{C}$.
For $r = \sum a_i\otimes b_i\in\mathfrak{g}\otimes\mathfrak{g}$, introduce
classical Yang---Baxter equation (CYBE) due to A.A. Belavin, V.G.~Drinfel'd~\cite{BelaDrin82} as
\begin{equation}\label{CYBE}
[r_{12},r_{13}]+[r_{12},r_{23}]+[r_{13},r_{23}] = 0,
\end{equation}
where
$$
r_{12} = \sum a_i\otimes b_i\otimes 1,\quad
r_{13} = \sum a_i\otimes 1\otimes b_i,\quad
r_{23} = \sum 1\otimes a_i\otimes b_i
$$
are elements from $\mathfrak{g}^{\otimes 3}$.

A tensor $r = \sum a_i\otimes b_i\in \mathfrak{g}\otimes\mathfrak{g}$
is called skew-symmetric if $\sum a_i\otimes b_i = - \sum b_i\otimes a_i$.
In \cite{BelaDrin82,Semenov83} it was actually shown that
given a skew-symmetric solution of CYBE on $\mathfrak{g}$,
a linear map $R\colon\mathfrak{g}\to\mathfrak{g}$ defined as
$R(x) = \sum \langle a_i,x\rangle b_i$
is an RB-operator of zero weight on~$\mathfrak{g}$.
Here $\langle \cdot,\cdot\rangle$ denotes the Killing form on $\mathfrak{g}$.

{\bf Example 4} \cite{Stolin}.
Up to conjugation and scalar multiple
unique skew-symmetric solution of CYBE on $\mathrm{sl}_2(\mathbb{C})$ is
$e\otimes h - h\otimes e$. Corresponding RB-operator is the following:
$R(e) = 0$, $R(f) = 4h$, $R(h) = -8e$.

Let $A$ be an associative algebra, a tensor
$r = \sum a_i\otimes b_i\in A\otimes A$ is a solution of
associative Yang---Baxter equation (AYBE, \cite{Zhelyabin,Aguiar00-2,Aguiar01}) if
\begin{equation}\label{AYBE}
r_{13}r_{12}-r_{12}r_{23}+r_{23}r_{13} = 0,
\end{equation}
where the definition of $r_{12},r_{13},r_{23}$ is the same as for CYBE.

{\bf Statement} \cite{Aguiar00}.
Let $r = \sum a_i\otimes b_i$ be a solution of AYBE on an associative algebra~$A$.
A linear map $R\colon A\to A$ defined as
$R(x) = \sum a_i x b_i$ is an RB-operator on $A$ of zero weight.

{\bf Example 5} \cite{Aguiar00-2}.
Up to conjugation, transpose and scalar multiple
all nonzero solutions of AYBE on $M_2(\mathbb{C})$ are
$(e_{11}+e_{22})\otimes e_{12}$;\
$e_{12}\otimes e_{12}$;\
$e_{22}\otimes e_{12}$;\
$e_{11}\otimes e_{12} - e_{12}\otimes e_{11}$.

In 1983 \cite{Semenov83}, M.A. Semenov-Tyan-Shansky
introduced modified Yang---Baxter equation (MYBE).
Let $L$ be a Lie algebra, $R$ be a linear map on $L$, then MYBE is
\begin{equation}\label{MYBE}
R(x)R(y) - R(R(x)y+xR(y)) = -xy.
\end{equation}

It is easy to check that $R$ is a solution of MYBE if and only if
$(R+\id)$ is an RB-operator on $L$ of the weight $-2$.
So, there is up to scalar multiple one-to-one correspondence
between the set of solutions of MYBE and RB-operators of nonzero weight.

\subsection{PBW-pair of varieties}

In \cite{PBW} it was introduced the notion of a PBW-pair
which generalizes the relation between associative and Lie algebras
given by famous Poincar\'{e}---Birkhoff---Witt theorem.

Let $\mathcal{V}$, $\mathcal{W}$ be varieties of algebras and
$\psi\colon \mathcal{V}\to \mathcal{W}$
be a such functor that maps $A\in \mathcal{V}$ to $\psi(A)\in \mathcal{W}$
preserving $A$ as vector space but changing the operations on $A$.
Universal enveloping algebra $U(A)$ is an image of $A$
of left adjoint functor to the $\psi$.
Defining on $U(A)$ a natural ascending filtration,
we get associated graded algebra $\mathrm{gr}\,U(A)$.

A pair of varieties $(\mathcal{V},\mathcal{W})$
with the functor $\psi\colon \mathcal{V}\to \mathcal{W}$
is called PBW-pair if $\mathrm{gr}\,U(A)\cong U(\mathrm{Ab}\,A)$.
Here $\mathrm{Ab}\,A$ denotes the vector space $A$
with trivial multiplicative operations.

\subsection{Free commutative RB-algebra}

Consider the one of possible constructions of
free commutative RB-algebra from \cite{GuoKeigher},
which is based on shuffle algebra.

Let $A$ be a commutative algebra over the field $\Bbbk$.
Denote by $\shuffle^+(A)$ the vector space
$\mathop{\oplus}\limits_{n=0}^\infty A^{\otimes n}$.
Define the bilinear operation $\diamond$ on $\shuffle^+(A)$
as follows: $a\diamond b$ for
$a = a_1\otimes \ldots \otimes a_m\in A^{\otimes m}$,
$b = b_1\otimes \ldots \otimes b_n\in A^{\otimes n}$
equals a sum of all tensors of the length $m+n$
which tensor factors are exactly
$a_i,b_j$, $i=1,\ldots,m$, $j=1,\ldots,n$,
moreover, the natural orders of $a_i$ and $b_j$ are preserved.

Give an explicit definition of the product $a\diamond b$
by induction on $m + n$. If at least one of the numbers $m,n$
equals 0, e.g. $n = 0$, then define $a\diamond b$
as a multiplication of the scalar $b\in \Bbbk$ on the tensor $a$ .
For $m + n = 2$, define $a\diamond b = a\otimes b + b\otimes a$.
For $m + n > 2$,
\begin{equation}\label{shuffl}
a\diamond b = \begin{cases}
a_1\otimes b + b_1\otimes (a_1\diamond b), & m=1,n\geq2, \\
a_1\otimes ((a_2\otimes\ldots \otimes a_m)\diamond b_1)
+ b_1\otimes a, & m\geq2,n=1, \\
a_1\otimes ((a_2\otimes\ldots \otimes a_m)\diamond b)
+ b_1\otimes (a\diamond (b_2\otimes\ldots\otimes b_n)), & m,n\geq2.
\end{cases}
\end{equation}

{\bf Example 6}.
Calculating $(a_1\otimes a_2)\diamond (b_1\otimes b_2)$,
we have $C_4^2 = 6$ summands:
\begin{multline*}
(a_1\otimes a_2)\diamond (b_1\otimes b_2)
 = a_1 \otimes (a_2\diamond (b_1\otimes b_2))
  + b_1\otimes ((a_1\otimes a_2)\diamond b_2) \\
 = a_1 \otimes a_2\otimes b_1\otimes b_2
 + a_1 \otimes b_1\otimes a_2\otimes b_2
 + a_1 \otimes b_1\otimes b_2\otimes a_2 \\
 + b_1\otimes a_1\otimes a_2\otimes b_2
 +  b_1\otimes a_1\otimes b_2\otimes a_2
 + b_1\otimes b_2\otimes a_1\otimes a_2.
\end{multline*}

Consider a tensor product $\shuffle(A) = A\otimes\shuffle^+(A)$
of algebras $\langle A,\cdot\rangle$ and $\langle \shuffle^+(A),\diamond\rangle$.
Define a linear map $R$ on $\shuffle(A)$ as follows
$$
R(a_0\otimes a) = \begin{cases}
1\otimes a_0\otimes a, & a\in A^{\otimes n}, n\geq1, \\
1\otimes a_0a, & a\in \Bbbk.
\end{cases}
$$

In \cite{GuoKeigher}, it was proved that
a subalgebra of $\shuffle(A)$ generated by $A\otimes 1$
is isomorphic to $\RB\Com(A)$, a free commutative RB-algebra
generated by $A$.
A free commutative RB-algebra generated by a set $X$
could be constructed as $\shuffle(\Bbbk[X])$ \cite{GuoKeigher}.
We will denote a free commutative RB-algebra of weight $\lambda$
generated by a set $X$ as $\RB_\lambda\Var\langle X\rangle$.
For short, an algebra $\RB_0\Var\langle X\rangle$
will be denoted as $\RB\Var\langle X\rangle$.

In \cite{GuoKeigher}, quasi-shuffle algebra
was used for the construction of free commutative RB-algebra
of nonzero weight.

In \cite{Cartier72}, P. Cartier actually proved that a linear base $T$
of $\RB_\lambda\Var\langle X\rangle$ could be construc\-ted by induction:

\noindent--- all monomials from $\Bbbk[X]$ lie in $T$,

\noindent--- if $u\in T$, then $R(u),R(u)w\in T$
for a monomial $w$ from $\Bbbk[X]$.

We will use the linear base of $\RB_\lambda\Com\langle X\rangle$ of P. Cartier \cite{Cartier72}
and will refer to it as a standard base.
Given a word $u$ from the standard base,
the number of appearances of the symbol $R$ in the notation of $u$
is called $R$-degree of the word $u$, denotation: $\deg_R(u)$.
Also define a degree of $u$ as follows:
$\deg(u) = \begin{cases}
 n, & u = x_{i_1}x_{i_2}\ldots x_{i_n},\ x_{i_j}\in X, \\
 1, & u = R(p), \\
 n+1, & u = R(p)x_{i_1}\ldots x_{i_n},\ x_{i_j}\in X.
\end{cases}$

\subsection{Precommutative algebra}

Precommutative algebra is an algebra which product satisfies the identity
\begin{equation}\label{preCom}
(x_1\succ x_2 + x_2\succ x_1)\succ x_3 = x_1\succ (x_2\succ x_3).
\end{equation}

A free precommutative algebra generated by a linear space $V$
could be constructed \cite{Dialg99} as reduced tensor algebra
$\bar{T}(V) = \oplus_{n\geq1} V^{\otimes n}$ with the product
\begin{multline}
v_1\otimes v_2\otimes\ldots \otimes v_p
\succ v_{p+1}\otimes \ldots \otimes v_{p+q} \\
 = ((v_1\otimes v_2\otimes\ldots \otimes v_p)
 \diamond (v_{p+1}\otimes \ldots \otimes v_{p+q-1}))\otimes v_{p+q}
\end{multline}
where $\diamond$ is defined by \eqref{shuffl}.

\subsection{Postcommutative algebra}

A postcommutative algebra is a linear space endowed with two bilinear products
$\succ$, $\perp$ such that $\perp$ is associative and commutative
and the following identities are fulfilled
\begin{equation}\label{PostCom}
\begin{gathered}
(x \succ y + y\succ x + x\perp y) \succ z = x \succ (y \succ z), \\
x\succ (y\perp z) = (x\succ y) \perp z,\quad (x \succ y)\perp z = y\perp (x \succ z).
\end{gathered}
\end{equation}

A free postcommutative algebra could be constructed
with the help of quasi-shuffle algebra \cite{Loday2007}.

\subsection{Embedding of Loday algebras into RB-algebras}

The common definition of the varieties of pre- and post-$\Var$-algebra
for a variety $\Var$ could be found in \cite{BBGN2012,GubKol2014}.

Given a commutative algebra $B$ with an RB-operator $R$ of zero weight,
the space $B$ with respect to the operation
\begin{equation}\label{eq:RB-PreOperations}
x\succ y = R(x)y
\end{equation}
is a precommutative algebra.

Given a commutative algebra $\langle B,\cdot\rangle$
with an RB-operator $R$ of unit weight,
we have that $\langle B,\succ,\cdot\rangle$
is a postcommutative algebra,
where the product $\succ$ is defined by \eqref{eq:RB-PreOperations}.
The case of RB-operator $R$ of any nonzero weight $\lambda$
is reduced to unit weight case as follows:
the map $\frac{1}{\lambda}R$ is an RB-operator of unit weight.

Given a precommutative algebra $\langle C,\succ\rangle$,
universal enveloping commutative RB-algebra $U$ of $C$
is an universal algebra in the class of all commutative RB-algebras
of zero weight such that there exists injective homomorphism from $C$ to $U^{(R)}$.
Analogously, we define an universal enveloping commutative RB-algebra
of a postcommutative algebra. The common denotation
of universal enveloping algebra: $U_{\RB}(C)$.

{\bf Theorem 1} \cite{GubKol2013}.
a) Any pre-$\Var$-algebra could be embedded into its universal envelo\-ping
RB-algebra of the variety $\Var$ and zero weight.

b) Any post-$\Var$-algebra could be embedded into its universal enveloping
RB-algebra of the variety $\Var$ and nonzero weight.

After Theorem 1, we have the natural question:
What does a linear base of universal enveloping RB-algebra
of a pre- or post-$\Var$-algebra look like for a variety $\Var$?
In the case of associative pre- and postalgebras, the question
appeared in \cite{Guo2011} and was solved in \cite{Gub2017}.
The current article is devoted to answer the question
in the commutative case.

In the article, the common method to construct universal enveloping is the following.
Let $X$ be a linear base of a precommutative algebra $C$.
We find a base of universal enveloping $U_{\RB}(C)$
as the special subset $E$ of the standard base of $\RB\Com\langle X\rangle$
closed under the action of RB-operator.
By induction, we define a product $*$
on the linear span of $E$ and prove its associativity.
Finally, we state universality of the algebra $\Bbbk E$.

In the case of postcommutative algebras,
as it was mentioned above, we will consider
universal enveloping commutative RB-algebra of unit weight.

\section{Universal enveloping Rota---Baxter algebra of pre\-commutative algebra}

Within the paragraph, we will construct universal enveloping RB-algebra
of an arbitrary precommutative algebra $\langle C,\succ\rangle$.

Let $B$ be a linear base of $C$.
Let us consider the algebra $A = \Bbbk[B]/I$ for
an ideal $I$ in $\Bbbk[B]$ generated by the set
$B' = \{(b\succ a)c-(b\succ c)a,\,a,b,c\in B\}$.
Here the expressions $b\succ a$ and $b\succ c$, $a,b,c\in B$,
equal the results of the products in $C$.
As $B$ is the linear base of $C$, the expressions are linear
combinations of the elements of $B$.
Denote by $\cdot$ the product in $A$.

Let us denote by $C(B)$ the set of all monomials in $\Bbbk[B]$.
Due to Gr\"{o}bner theory,
there exists such set $E_0\subset C(B)$ that $\bar{E}_0$ ---
the image of $E_0$ under the factorization of $\Bbbk[B]$
by the ideal $I$ --- is a base of $A$; moreover,
for any decomposition of an element $v\in E_0$ into a concatenation
$v = v_1 v_2$ of nonempty $v_1,v_2$ we have $v_1,v_2\in E_0$.

An element $R(u)$ of the standard base is called good if $u\neq b,R(x)b$ for $b\in B$.

Let us define by induction $Envelope$-words (shortly $E$-words),
a subset of the standard base of $\RB\Com\langle B\rangle$:

1) elements of $E_0$ are $E$-words of the type 1;

2) given an $E$-word $u$,
we define $R(u)$ as an $E$-word of the type 2;

3) given an $E$-word $x$, we define $R(x)w$ for $w\in E_0$
as an $E$-word of the type 3 if $R(x)$ is good.

{\bf Theorem 2}.
The set of all $E$-words forms a linear base of
universal enveloping commutative RB-algebra of $C$.

{\bf Lemma 1}.
Let $D$ denote a linear span of all $E$-words.
One can define such bilinear commutative operation $*$ on the space $D$ that
($k$--$l$ denotes below the condition on the product $v*u$,
where $v$ is an $E$-word of the type $k$ and $u$ is an $E$-word of the type $l$.)

1--1: given $v,u\in E_0$, we have
\begin{equation}\label{def:1--1}
v*u = v\cdot u.
\end{equation}

1--2: given $v = w'a\in E_0$, $a\in B$, $w'\in E_0\cup \{\emptyset\}$,
an $E$-word $u = R(p)$ of the type 2, we have
\begin{equation}\label{def:1--2}
v*u
 = \begin{cases}
  w'\cdot (b\succ a), & u = R(b), b\in B, \\
  R(x)(w'\cdot(b\succ a))-R(x*R(b))w'a, & u = R(R(x)b), b\in B, \\
  R(p)w'a, & else.
    \end{cases}
\end{equation}
(If $w'=\emptyset$, by $w'\cdot(b\succ a)$ we mean $b\succ a$.)

1--3: given $v = w_1\in E_0$, an $E$-word $u = R(x)w_2$ of the type 3, $w_2\in E_0$, we have
\begin{equation}\label{def:1--3}
v*u = R(x)(w_1\cdot w_2).
\end{equation}

2--2: given $E$-words $v = R(p)$, $u = R(s)$ of the type 2, we have
\begin{equation}\label{def:2--2}
v*u = R( R(p)*s + R(s)*p ).
\end{equation}

2--3: given an $E$-word $v = R(x)$ of the type 2 and $E$-word
$v = R(y)w$, $w = w'a\in E_0$, $w'\in E_0\cup\emptyset$, $a\in B$,
of the type 3, we have
\begin{equation}\label{def:2--3}
v*u =
  \begin{cases}
   R(y)(w'\cdot(b\succ a)), & x=b\in B, \\
   (R(y)*R(z))(w'\cdot(b\succ a)) - (R(y)*R(z*R(b)))w & x=R(z)b,\,b\in B, \\
   (R(y)*R(x))w, & else.
   \end{cases}
\end{equation}

3--3: given $E$-words $v = R(x)w_1$, $u = R(y)w_2$,
$w_1,w_2\in E_0$, of the type 3, we have
\begin{equation}\label{def:3--3}
v*u = (R(x)*R(y))(w_1\cdot w_2).
\end{equation}

{\sc Proof}.
Let us define the operation $*$ with the prescribed conditions
for $E$-words $v,u$ by induction on $r = \deg_R(v) + \deg_R(u)$.
For $r = 0$ define $v*u = v\cdot u$, $v,u\in E_0$,
which satisfies the condition 1--1.

For $r = 1$, let us define $v*u$ as follows:

1--2: given $v = w_1 = w_1'a\in E_0$, $a\in B$,
an $E$-word $u = R(w_2)$ of the type 2, $w_2\in E_0$,
$$
v*u = \begin{cases}
  w_1'\cdot(b\succ a), & w_2 = b \in B, \\
  R(w_2)w_1, & else.
               \end{cases}
$$

1--3: given $v = w_1\in E_0$, an $E$-word $u = R(w_2)w_3$
of the type 3, $w_2\in E_0\setminus B$, $w_3\in E_0$,
$v*u = R(w_2)(w_1\cdot w_3)$.

For $r>1$, let us define $v*u$ as follows:

1--2: given $v = w'a\in E_0$, $a\in B$, $w'\in E_0\cup \{\emptyset\}$,
an $E$-word $u = R(p)$ of the type 2,
$$
v*u = \begin{cases}
  R(x)(w'\cdot(b\succ a))-R(x*R(b))w'a, & u = R(R(x)b), b\in B, \\
  R(p)w'a, & else.
               \end{cases}
$$
It is correct to write $R(x*R(b))w'a$. Indeed,
as $p = R(x)b$ is $E$-word, so $R(x)$ is good.
We have three variants: $x = w\in E_0$, $\deg (w)\geq2$,
$x = R(p)$, or $x = R(p)w'$, $w'\in E_0$. In all cases
$x*R(b)$ is a sum of $E$-words of the form
different from $c,R(s)c$ for $c\in B$.

1--3, 2--2, 2--3, 3--3: by the formulas \eqref{def:1--3}, \eqref{def:2--2},
\eqref{def:2--3}, \eqref{def:3--3} respectively.

Prove that the definition of the product by \eqref{def:2--3} is correct.
By the reasons stated above $R(z*R(b))$ is a sum of good $E$-words of the type 2.
Thus, it is remains to check that $R(p)*R(s)$ is a sum of good $E$-words of the type 2
provided $R(p)$ and $R(s)$ are so. Actually it is enough to verify that
$R(R(p)*s)$ is a sum of good $E$-words of the type 2.
It is easy to state the last fact by consideration of three variants of $s$.

We put the product $u*v$ in the cases 2--1, 3--1, 3--2 equal to $v*u$.

{\bf Lemma 2}.
The space $D$ with the operations $*$, $R$ is an RB-algebra.

{\sc Proof}.
It follows from \eqref{def:2--2}.

{\bf Lemma 3}.
The relation $R(b)*a = b\succ a$ holds in $D$ for every $a,b\in B$.

{\sc Proof}.
It follows from Lemma 1, the first case of \eqref{def:1--2}.

{\bf Lemma 4}.
The operation $*$ on $D$ is commutative and associative.

{\sc Proof}.
The commutativity follows from Lemma 1.

Given $E$-words $x,y,z$, let us prove associativity
$$
(x,y,z) = (x*y)*z - x*(y*z) = 0
$$
by inductions on two parameters:
at first, on summary $R$-degree $r$ of the triple $x,y,z$,
at second, on summary degree $d$ of $x,y,z$.

For $r = 0$, we have $(x,y,z) = (x\cdot y)\cdot z - x\cdot (y\cdot z) = 0$, $x,y,z\in E_0$,
as the product $\cdot$ is associative in the algebra $A$.

Let $r>0$ and suppose that associativity
for all triples of $E$-words with the less summary $R$-degree is proven.
Consider $d = 3$.

2--1--1: let $x = R(p)$, $y = a$, $z = c$, $a,c\in B$.

a) Let $x = R(b)$, $b\in B$.
Remember that $A = \Bbbk[B]/I$ for
$I = \langle (d\succ a)c-(d\succ c)a,\,a,d,c\in B\rangle$.
Define a map $\vdash\colon \Bbbk[B]\to A$
on monomials as follows:
\begin{equation}\label{eq:vdashv}
\vdash (aw) = (b\succ a)\cdot w,
\end{equation}
where $b\succ a$ is a product in $C$ and equals to a linear combination of elements of $B$.
The map $\vdash$ is well-defined because of the view of $I$.

Let us check that $\vdash (I) = 0$. By \eqref{preCom} and the definition of $I$,
\begin{multline*}
b*((d\succ a)cu-(d\succ c)au) \\
 = (b\succ(d\succ a))\cdot cu
  - (b\succ(d\succ c))\cdot au \\
 = u\cdot((b\succ(d\succ a))\cdot c -(b\succ(d\succ c))\cdot a) \\
 = u\cdot(((b\succ d + d\succ b)\succ a)\cdot c
         - ((b\succ d + d\succ b)\succ c)\cdot a)
 = 0.
\end{multline*}

Thus, the map $\vdash$ could be considered as map $\vdash\colon A\to A$
and it coincides on $E_0$ with the map $z\to R(b)*z$. Finally,
$$
(x,y,z) = (R(b),a,c)
  = (b\succ a)\cdot c - R(b)*(a\cdot c)
  = (b\succ a)\cdot c - b\vdash ac
  = 0.
$$

b) The case $x = R(R(y)b)$, $b\in B$, could be reduced from a).

c) If $x$ is good, then associativity follows from Lemma 1.

1--2--1: let $x = a$, $y = R(u)$, $z = c$, $a,c\in B$.
If $y$ is good, the associativity follows by \eqref{def:1--2} and \eqref{def:1--3}.

Let $y = R(b)$, $b\in B$. We have
$(x,y,z) = (b\succ a)\cdot c - (b\succ c)\cdot a = 0$, as
we deal with the product in $A$.

Let $y = R(R(u)b)$, $b\in B$. From one hand,
\begin{multline*}
(x*y)*z
 = (a*R(R(u)b))*c \\
 = (R(u)(b\succ a)-R(u*R(b))a)*c \\
 = R(u)((b\succ a)\cdot c) - R(u*R(b))(a\cdot c).
\end{multline*}
From another hand,
\begin{multline*}
x*(y*z)
 = a*(R(R(u)b))*c) \\
 = a*(R(u)(b\succ c)- R(u*R(b))c) \\
 = R(u)((b\succ c)\cdot a) - R(u*R(b))(a\cdot c).
\end{multline*}
Associativity follows directly.

2--2--1: let $x = R(u)$, $y = R(v)$, $z = a\in B$.
If $y$ is good, then associativity follows from \eqref{def:1--2} and \eqref{def:2--2}.

a) Let $y = R(c)$, $c\in B$.
a1) If $x = R(b)$, $b\in B$, then by \eqref{preCom}
\begin{multline*}
(x,y,z)
 = R(R(b)*c + b*R(c))*a - R(b)*(c\succ a) \\
 = R(b\succ c + c\succ b)*a - (b\succ (c\succ a)) \\
 = ((b\succ c + c\succ b)\succ a) - (b\succ (c\succ a)) = 0.
\end{multline*}

a2) If $x = R(R(p)b)$, $b\in B$, then compute
\begin{multline}\label{Com:2-2-1a2}
(x,y,z)
 = (R(R(p)b)*R(c))*a - R(R(p)b)*(c\succ a) \\
 = R( R(R(p)b)*c )*a + R(R(p)b*R(c))*a
   + R(p*R(b))*(c\succ a) - R(p)(b\succ(c\succ a)) \\
 = R(R(p)(b\succ c))*a - R(R(p*R(b))c)*a + R(R(R(p)c)*b)*a \\
   + R(R(p*R(c))b)*a + R(p*R(b))*(c\succ a) - R(p)(b\succ(c\succ a)) \\
 =    R(p)((b\succ c)\succ a) - R(p*R((b\succ c)))a \\
 +R((p*R(b))*R(c))a - R(p*R(b))(c\succ a) \\
 +R(R(p)(c\succ b))*a - R(R(p*R(c))b)*a \\
 +R(p*R(c))(b\succ a) - R((p*R(c))*R(b))a \\
 +R(p*R(b))*(c\succ a) - R(p)(b\succ(c\succ a)).
\end{multline}
Substituting the equalities
\begin{gather*}
R(R(p)(c\succ b))*a
 = R(p)((c\succ b)\succ a) - R(p*R(c\succ b))a, \\
R(R(p*R(c))b)*a
 = R(p*R(c))(b\succ a) - R((p*R(c))*R(b))a,
\end{gather*}
into \eqref{Com:2-2-1a2} and applying \eqref{preCom}
and the equality $(p,R(b),R(c)) = 0$ holding by induction on $r$,
we obtain zero.

a3) Let $x$ be good, then
\begin{equation}\label{Com:2-2-1a3}
(x,y,z) = R(R(u)c + u*R(c))*a - R(u)(c\succ a).
\end{equation}

Calculating $R(R(u)c)*a$ in \eqref{Com:2-2-1a3} by \eqref{def:1--2},
we obtain $(x,y,z) = 0$.

b) Let $y = R(R(t)c)$, $c\in B$.
\begin{multline}\label{Com:2-2-1b}
(x,y,z)
 = (R(u)*R(R(t)c))*a - R(u)*(R(R(t)c)*a) \\
 = R(R(u)*R(t)c)*a + R(u*R(R(t)c))*a
 - R(u)*(R(t)(c\succ a)-R(t*R(c))a) \\
 = R(R(R(u)*t)c)*a + R(R(u*R(t))*c)*a
 + R(u*R(R(t)c))*a + R(R(u)*(t*R(c)))*a \\
 + R(u*R(t*R(c)))*a - R(R(u)*t)(c\succ a) - R(u*R(t))*(c\succ a).
\end{multline}
The correctness of the notation $R(R(u)*t)e)$, $e\in B$, easily follows from the
fact that $R(t)$ is good. Write down the first summand of the RHS of \eqref{Com:2-2-1b}
$$
R(R(R(u)*t)c)*a = R(R(u)*t)(c\succ a) - R((R(u)*t)*R(c))*a
$$
and substitute it into \eqref{Com:2-2-1b}.
Applying the equality $(R(u),t,R(c)) = 0$ holding by induction on $r$,
we have by the case a)
\begin{multline*}
(x,y,z) = R(R(u*R(t))*c)*a + R(u*R(t)*R(c))*a - R(u*R(t))*(c\succ a) \\
 = (R(u*R(t)),R(c),a) = 0.
\end{multline*}

2--1--2: associativity follows from the case 2--2--1.

2--2--2: by induction we have
\begin{equation}\label{Com:2-2-2}
(R(x),R(y),R(z)) = R( (R(x),R(y),z) + (R(x),y,R(z)) + (x,R(y),R(z)) = 0.
\end{equation}

Let $d > 3$. The proof for the cases
2--1--1, 1--2--1, 2--2--1, 2--1--2, 2--2--2 is analogous to the written above.

3--2--1: let $x = R(u)w_1$, $y = R(t)$, $z = w_2 = aw_2'\in E_0$, $a\in B$.
If $t = b\in B$, then by induction on $d$ we have
\begin{gather*}
(x*y)*z = ((R(u)*R(b))*w_1)*w_2
 = ((R(u)*R(b))*w_2)*w_1, \\
x*(y*z) = R(u)(w_1\cdot(b\succ a)\cdot w_2')
 = (R(u)*(R(b)*w_2))*w_1.
\end{gather*}
Therefore, $(x,y,z) = (R(u),R(b),w_2)*w_1 = 0$ by induction on $d$.

Let $y = R(R(v)b)$, $b\in B$.
\begin{equation}\label{Com:3-2-1L}
(x*y)*z = ((R(u)*(R(R(v)b)))*w_1)*w_2,
\end{equation}
\begin{multline}\label{Com:3-2-1R}
x*(y*z) =
 R(u)w_1*(R(R(v)b)*w_2) \\
 = R(u)w_1*(R(v)((b\succ a)\cdot w_2') -R(v*R(b))w_2 ) \\
 = (R(u)*R(v))((b\succ a)\cdot w_2'\cdot w_1) -(R(u)*R(v*R(b))(w_1\cdot w_2).
\end{multline}
Applying induction on $d$ for the first summand
of the RHS of \eqref{Com:3-2-1R}, rewrite
\begin{multline*}
(R(u)*R(v))((b\succ a)\cdot w_2'\cdot w_1)
 = (R(u)*R(v))*((R(b)*w_2)*w_1) \\
 = ((R(u)*R(v))*(R(b)*w_2))*w_1
 = ((R(u)*R(v))*R(b))*(w_1\cdot w_2).
\end{multline*}
Analogously the RHS of \eqref{Com:3-2-1L} equals
$$
(( R(u)*(R(R(v)b)))*w_1)*w_2
 = (R(u)*(R(R(v)b))*(w_1\cdot w_2).
$$
Finally, we have
\begin{multline*}
(x*y)*z - x*(y*z) \\
 = (R(u)*(R(R(v)b))*(w_1\cdot w_2)
 + (R(u)*R(v*R(b))(w_1\cdot w_2)
 - ((R(u)*R(v))*R(b))*(w_1\cdot w_2) \\
 = -(R(u),R(v),R(b))*(w_1\cdot w_2) = 0.
\end{multline*}

If $y = R(t)$ is good, then
\begin{multline*}
(x*y)*z = ((R(u)*R(t))w_1)*w_2
 = (R(u)*R(t))w_1w_2 \\
 = (R(u)w_1)*(R(t)w_2) = x*(y*z).
\end{multline*}

3--1--2: the proof is analogous to the one in the case 2--1--2.

2--3--1: associativity follows from the cases 2--2--1 and 2--1--1.

3--1--1, 1--3--1, 3--3--1, 3--1--3:
associativity follows from \eqref{def:1--1}, \eqref{def:1--3}.

3--3--2: let $x = R(u_1)w_1$, $y = R(u_2)w_2$, $z = R(u_3)$.
If $z$ is not good, then the proof follows from the proofs
in the cases 1--1--2 and 2--2--2. If $z$ is good, from one hand,
$$
(x*y)*z = ((R(u_1)*R(u_2))w_1w_2)*R(u_3)
= ((R(u_1)*R(u_2))*R(u_3))*w_1w_2.
$$
From another hand, applying induction on $d$, we derive
\begin{multline*}
x*(y*z) = (R(u_1)w_1)*((R(u_2)*R(u_3))*w_2) \\
 = (R(u_1)w_1*(R(u_2)*R(u_3)))*w_2
 = (R(u_1)*(R(u_2)*R(u_3)))*w_1w_2.
\end{multline*}
Thus, associativity follows from the case 2--2--2.

3--2--2, 3--2--3 and 2--3--2: associativity
follows by the same inductive reasons as in the case 3--3--2.

3--3--3: associativity follows from \eqref{def:3--3}, \eqref{Com:2-2-2},
and associativity of the product in $A$.

{\sc Proof of Theorem 2}.
Let us prove that the algebra $D$ is exactly universal envelo\-ping algebra
for the precommutative algebra $C$, i.e., is isomorphic to the algebra
$$
U_{\RB\Com}(C) = \RB\Com\langle B|b\succ a = R(b)a,\,a,b\in B\rangle.
$$
By the construction, the algebra $D$ is generated by $B$.
Thus, $D$ is a homomorphic image of a homomorphism $\varphi$ from $U_{\RB\Com}(C)$.
Prove that all basic elements of $U_{\RB\Com}(C)$ are linearly expressed by $D$,
this leads to nullity of kernel of $\varphi$ and $D\cong U_{\RB\Com}(C)$.

As the equality
$$
(b\succ a)c = R(b)ac = R(b)ca = (b\succ c)a
$$
is satisfied on $U_{\RB\Com}(C)$,
the space $U_{\RB\Com}(C)$ is a subspace of
$\RB\Com\langle B|(b\succ a)c = (b\succ c)a,\,a,b,c\in B\rangle$.
It is well-known \cite{Guo2011} that algebras
$\RB\Com\langle B\rangle$ and $\RB\Com(\Bbbk[B])$ coincide.
Consider the map
$\varphi\colon B\to \RB\Com(\Bbbk[B])/J$,
where $J$ is the ideal generated by the set
$B' = \{(b\succ a)c = (b\succ c)a,\,a,b,c\in B\}$
and closed under RB-operator,
the map $\psi$ equal the composition of trivial maps
$\psi_1\colon B\to A$ and $\psi_2\colon A\to\RB\Com(A)$,
where as above,
$$
A = \Bbbk[B]/I,\ I =\langle B'\rangle_{\mathrm{id}}\triangleleft\Bbbk[B],\quad
\RB\Com(A) = \RB\Com\langle E_0|w_1 w_2 = w_1\cdot w_2\rangle.
$$
The base $S$ of $\RB\Com(A)$ \cite{Guo2011}
could be constructed by induction:
at first, $E_0\subset S$;
at second, if $u\in S$, then
$R^k(R(u)w)\in S$ for $w\in E_0$, $k\geq0$.

From $\ker \psi\subseteq \ker \varphi$,
we have the injective embedding $U_{\RB\Com}(C)$ into $\RB\Com(A)$.

It remains to show by \eqref{eq:RB-PreOperations}
that the complement $E'$ of the set of all $E$-words
in the standard base of $\RB\Com(A)$ is linearly expressed via $E$-words.
Applying the inductions on the $R$-degree and
the degree of base words in $\RB\Com(A)$, the relations
\begin{gather*}
R(a)u = \begin{cases}
u'(a\succ b), & u = bu', b\in B, \\
R(R(a)t)u' + R(aR(t))u', & u = R(t)u'; \end{cases} \\
R(R(u)b)a = R(u)R(b)a - R(R(b)u)a
 = R(u)(b\succ a) - R(R(b)u)a,\quad a,b\in B,
\end{gather*}
we prove that the elements of $E'$
are linearly expressed via $E$-words.
The theorem is proved.

Given a precommutative algebra $C$,
$U_0(C)$ denotes a linear span of all $E$-words
of zero $R$-degree in $U_{\RB}(C)$.

{\bf Example 7}.
Let an algebra $A$ be a direct sum of $n$ copies of $\Bbbk$.
Consider $A$ as a precommutative algebra under the operations
$a\succ b = ab$, $a\prec b = 0$.
From $A^2 = A$ we conclude $U_0(A) = A$.

{\bf Example 8}.
Let $B$ be a $n$-dimensional vector space over the field $\Bbbk$
with trivial operations $a\succ b = a\prec b = 0$.
So, $U_0(B) = \Bbbk[B]$.

{\bf Corollary 1}.
The pair of varieties $(\RB\Com, \pre\Com)$ is not a PBW-pair.

{\sc Proof}.
Universal enveloping commutative RB-algebra of
finite-dimensional precom\-mutative algebras $A$ and $B$ (Examples 7 and 8)
of the same dimension are not isomorphic, else the spaces
$U_0(A)$ and $U_0(B)$ were isomorphic as vector spaces. But we have
$$
\dim U_0(A) = \dim A = n < \dim U_0(B) = \dim \Bbbk[B] = \infty.
$$

Thus, the structure of universal enveloping commutative RB-algebra
of a precommu\-ta\-tive algebra $C$ depends on the operation $\succ$ of $C$.

\section{Universal enveloping Rota---Baxter algebra of post\-commutative algebra}

In the paragraph, we will construct universal enveloping RB-algebra
of an arbitrary postcommutative algebra $\langle C,\succ,\perp\rangle$.
Let $B$ be a linear base of $C$.

Let us define by induction $Envelope$-words (shortly $E$-words),
a subset of the standard base of $\RB\Com\langle B\rangle$:

1) elements of $B$ are $E$-words of the type 1;

2) given an $E$-word $u$, we define $R(u)$ as an $E$-word of the type 2;

3) given an $E$-word $u$, we define $R^2(x)a$ for $a\in B$
as an $E$-word of the type 3.

{\bf Theorem 3}.
The set of all $E$-words forms a linear base of
universal enveloping commutative RB-algebra of $C$.

{\bf Lemma 5}.
Let $D$ denote a linear span of all $E$-words.
One can define such bilinear commutative operation $*$ on the space $D$ that
($k$--$l$ denotes below the condition on the product $v*u$,
where $v$ is an $E$-word of the type $k$ and $u$ is an $E$-word of the type $l$.)

1--1: given $v = a$, $u = b$, $a,b\in B$, we have
\begin{equation}\label{def:1--1post}
v*u = a\perp b.
\end{equation}

1--2: given $v = a\in B$, an $E$-word $u = R(p)$ of the type 2, we have
\begin{equation}\label{def:1--2post}
v*u
 = \begin{cases}
  (b\succ a), & u = R(b), b\in B, \\
  R^2(x)(b\succ a)-R(R(x)*R(b))a-R(R(x)*b)*a, & u = R(R^2(x)b), b\in B, \\
  R(p)a, & else.
    \end{cases}
\end{equation}

1--3: given $v = a\in B$, an $E$-word $u = R^2(x)b$ of the type 3, $b\in B$, we have
\begin{equation}\label{def:1--3post}
v*u = R^2(x)(a\perp b).
\end{equation}

2--2: given $E$-words $v = R(p)$, $u = R(s)$ of the type 2, we have
\begin{equation}\label{def:2--2post}
v*u = R( R(p)*s + R(s)*p + s*p ).
\end{equation}

2--3: given an $E$-word $v = R^2(x)a$ of the type 3 and $E$-word
$u = R(y)$ of the type 2, $a\in B$, we have
\begin{equation}\label{def:2--3post}
v*u = R(R(x)*R(y))*a.
\end{equation}

3--3: given $E$-words $v = R^2(x)a$, $u = R^2(y)b$, $a,b\in B$, of the type 3, we have
\begin{equation}\label{def:3--3post}
v*u = (R(R(x))*R(R(y)))(a\perp b).
\end{equation}

{\sc Proof}.
Proof is analogous to the proof of Lemma 1.

{\bf Lemma 6}.
The algebra $D$ is an enveloping commutative
RB-algebra of unit weight for $C$.

{\sc Proof}.
In general, the proof is analogous to the proof of Lemmas 2--4.

The key moment is to prove an associativity.
We prove that
$(x,y,z) = (x*y)*z - x*(y*z) = 0$ for a triple $x,y,z$ of the types $k$--$l$--$m$.
We proceed by inductions on $r = \deg_R(x)+\deg_R(y)+\deg_R(z)$
and $d = \deg(x)+\deg(y)+\deg(z)$.
We will consider only the cases which proofs are technically different
from the proof of Lemma 4.

2--1--1.
Consider the case $x = R(b)$, $y = a$, $z = c$ for $a,b,c\in B$. Then by \eqref{PostCom}
$$
(x,y,z)
 = (R(b)*a)*c - R(b)*(a*c)
 = (b\succ a)\perp c - b\succ (a\perp c) = 0.
$$

2--2--1. a) Let $x = R(R^2(p)b)$ be an $E$-word of the type 2, $y = R(c)$, $z = a$; $a,b,c\in B$.
\begin{multline}\label{PostCom2-2-1a-R}
x*(y*z) = R( R^2(p)b)*(c\succ a) \\
= R^2(p)(b\succ(c\succ a)) - R(R(p)*R(b))*(c\succ a) - R(R(p)*b)*(c\succ a).
\end{multline}
\begin{equation}\label{PostCom2-2-1a-L}
(x*y)*z = R( R(R^2(p)b)*c )*a + R( R^2(p)b*R(c) )*a
 + R( R^2(p)(b\perp c))*a.
\end{equation}
Applying induction on $r$, write down the RHS of \eqref{PostCom2-2-1a-L} more detailed
\begin{multline}\label{PostCom2-2-1a-L1}
R( R(R^2(p)b)*c )*a
 = R( R^2(p)(b\succ c) - R(R(p)*R(b))*c - R(R(p)*b)*c )*a \\
 = R^2(p)((b\succ c)\succ a) - R(R(p)*R(b\succ c))*a - R(R(p)*(b\succ c))*a \\
 - R(R(p)*R(b))*(c\succ a) + R(R(p)*R(b)*R(c))*a + R(R(p)*(b\succ c))*a \\
 - R(R(p)*b)*(c\succ a) + R(R(p)*(c\succ b)) *a + R(R(p)*(b\perp c))*a \\
{=}  R(R(p)*R(b)*R(c))*a  - R(R(p)*R(b\succ c))*a
  + R(R(p)*(c\succ b)) *a + R(R(p)*(b\perp c))*a \\
  - R(R(p)*R(b))*(c\succ a) - R(R(p)*b)*(c\succ a) + R^2(p)((b\succ c)\succ a);
\end{multline}
\begin{multline}\label{PostCom2-2-1a-L2}
R( R^2(p)b*R(c) )*a
 = R(R(R^2(p)c)*b)*a + R(R(R(p)*R(c))*b)*a + R(R(R(p)*c)*b)*a \\
 = R(R^2(p)(c\succ b))*a - R(R(R(p)*R(c))*b)*a - R(R(R(p)*c)*b)*a \\
 + R(R(p)*R(c))*(b\succ a) - R(R(p)*R(c)*R(b))*a - R(R(p)*(c\succ b))*a \\
 + R(R(p)*c)*(b\succ a) - R(R(p)*(b\succ c))*a - R(R(p)*(b\perp c))*a;
\end{multline}
\begin{multline}\label{PostCom2-2-1a-L3}
R( R^2(p)(b\perp c))*a  \\
 = R^2(p)*((b\perp c)\succ a) - R(R(p)*R(b\perp c))*a - R(R(p)*(b\perp c))*a.
\end{multline}
Write down the second row of \eqref{PostCom2-2-1a-L2} as
\begin{multline}\label{PostCom2-2-1a-L2b}
R(R^2(p)(c\succ b))*a - R(R(R(p)*R(c))*b)*a - R(R(R(p)*c)*b)*a \\
 = R^2(p)((c\succ b)\succ a) - R(R(p)*R(c\succ b))*a - R(R(p)*(c\succ b))*a \\
 - R(R(p)*R(c))*(b\succ a) + R(R(p)*R(c)*R(b))*a + R(R(p)*(c\succ b))*a \\
 - R(R(p)*c)*(b\succ a) + R(R(p)*(b\succ c))*a + R(R(p)*(b\perp c))*a.
\end{multline}

Substituting \eqref{PostCom2-2-1a-L2b} into \eqref{PostCom2-2-1a-L2}, we obtain
\begin{multline}\label{PostCom2-2-1a-L2fin}
R( R^2(p)b*R(c) )*a \\
 = R^2(p)((c\succ b)\succ a) - R(R(p)*R(c\succ b))*a - R(R(p)*(c\succ b))*a.
\end{multline}

Subtracting \eqref{PostCom2-2-1a-R} from the sum of
\eqref{PostCom2-2-1a-L1}, \eqref{PostCom2-2-1a-L3}, \eqref{PostCom2-2-1a-L2fin},
we get zero by induction.

b) $x = R(u)$, $y = R(R^2(t)b)$, $z = a$; $a,b\in B$.
Applying induction on $r$, we have
\begin{multline}\label{PostCom2-2-1bL}
(x*y)*z
 = R(R(u)*R^2(t)b)*a + R(u*R(R^2(t)b))*a + R(u*R^2(t)b)*a  \\
 = R(R(u)*R(t) + u*R^2(t) + u*R(t))*(b\succ a) \\
 - R((R(u)*R(t) + u*R^2(t) + u*R(t))*R(b))*a \\
 - R((R(u)*R(t) + u*R^2(t) + u*R(t))*b)*a \\
 + R(u*R(R^2(t)b))*a + R(u*R^2(t)b)*a.
\end{multline}
\begin{multline}\label{PostCom2-2-1bR}
x*(y*z)
 = R(u)*(R^2(t)*(b\succ a) - R(R(t)*R(b))a - R(R(t)*b)*a) \\
 = R(R(u)*R(t) + u*R^2(t) + u*R(t))*(b\succ a) \\
 - R(R(u)*R(t)*R(b))*a - R(u*R(R(t)*R(b)))*a - R(u*R(t)*R(b))*a \\
 - R(R(u)*R(t)*b)*a - R(u*R(R(t)*b))*a - R(u*R(t)*b)*a.
\end{multline}
Subtracting \eqref{PostCom2-2-1bR} from \eqref{PostCom2-2-1bL},
we obtain zero by induction.

{\sc Proof of Theorem 3}.
The proof is analogous to the proof of Theorem 2.

{\bf Corollary 2}.
The pair of varieties $(\RB_\lambda\Com, \post\Com)$ is a PBW-pair.

It is well-known \cite{Hopf} that one can define a Hopf algebra structure
on a free commutative RB-algebra.

{\bf Question}.
Does there exist a Hopf algebra structure on universal enveloping
commu\-ta\-tive RB-algebra of a pre- and postcommutative algebra?

\noindent
Gubarev Vsevolod \\
Sobolev Institute of Mathematics of the SB RAS \\
Acad. Koptyug ave., 4 \\
Novosibirsk State University \\
Pirogova str., 2 \\
Novosibirsk, Russia, 630090\\
{\it E-mail: wsewolod89@gmail.com}

\end{document}